%%%% prevent double loading:
\expandafter\ifx\csname mthreemacsloaded\endcsname\relax\else \fi

\magnification1100
\input amstex

%%% Hack of Plain TeX correction and style macros 
%%% written by Walter Neumann and Larry Siebenmann:

 \catcode`\@=11
 \let\wlog@ld\wlog
 \def\wlog#1{\relax}

 \newif\ifIN@
 \def\m@rker{\m@@rker}
 \def\IN@{\expandafter\INN@\expandafter}
 \long\def\INN@0#1@#2@{\long\def\NI@##1#1##2##3\ENDNI@
    {\ifx\m@rker##2\IN@false\else\IN@true\fi}%
     \expandafter\NI@#2@@#1\m@rker\ENDNI@}
  \newtoks\Initialtoks@  \newtoks\Terminaltoks@
  \def\SPLIT@{\expandafter\SPLITT@\expandafter}
  \def\SPLITT@0#1@#2@{\def\TTILPS@##1#1##2@{%
     \Initialtoks@{##1}\Terminaltoks@{##2}}\expandafter\TTILPS@#2@}
  \newtoks\Trimtoks@

 \def\ForeTrim@{\expandafter\ForeTrim@@\expandafter}
 \def\ForePrim@0 #1@{\Trimtoks@{#1}}
 \def\ForeTrim@@0#1@{\IN@0\m@rker. @\m@rker.#1@%
     \ifIN@\ForePrim@0#1@%
     \else\Trimtoks@\expandafter{#1}\fi}
 
  \def\Trim@0#1@{%
      \ForeTrim@0#1@%
      \IN@0 @\the\Trimtoks@ @%
        \ifIN@
             \SPLIT@0 @\the\Trimtoks@ @\Trimtoks@\Initialtoks@
             \IN@0\the\Terminaltoks@ @ @%
                 \ifIN@
                 \else \Trimtoks@ {FigNameWithSpace}%
                 \fi
        \fi
      }

  %%% Math Bolds
  \font\titlebold=cmbx12 scaled 1200
  \font\twelvebold=cmbx12
  \font\tenbold=cmbx10
  \font\ninebold=cmbx9
  \font\sevenbold=cmbx7
  \font\fivebold=cmbx5

  \input amssym.def \input amssym
  %%% point sizes not loaded by amssym.def:
     \font\titlemsa=msam10 at 14.4pt
     \font\titlemsb=msbm10 at 14.4pt
     \font\titleeufm=eufm10 at 14.4pt
     \font\twelvemsa=msam10 scaled 1200
     \font\twelvemsb=msbm10 scaled 1200
     \font\twelveeufm=eufm10 scaled 1200
     \font\ninemsa=msam9
     \font\ninemsb=msbm9
     \font\nineeufm=eufm9

   %%% Cyrillic fonts (for accents and input, see ams cyr doc)
   \ifx\cyrfam\undefined
   \else
     \immediate\write16{}%
     \message{ !!! cyr fonts already defined. !!! }
     \message{ --- edit out superfluous font defs? }
   \fi
   \newfam\cyrfam
       \font\titlecyr=wncyr10 scaled 1440 %%% no caps?
       \font\twelvecyr=wncyr10 scaled 1200
       \font\tencyr=wncyr10
       \font\ninecyr=wncyr9
       \font\sevencyr=wncyr7
       \font\sixcyr=wncyr6

   %%% Euler script fonts (replacing caligraphic):
   \newfam\eusmfam
       \font\titleeusm=eusm10 scaled 1440
       \font\twelveeusm=eusm10 scaled 1200
       \font\teneusm=eusm10
       \font\nineeusm=eusm9
       \font\seveneusm=eusm7
       
       \font\fiveeusm=eusm5

\let\Cal\cal

 %%% Some fonts not loaded by plain
    \font\ninemrm=cmr9 %% new name for 9 pt math roman
    \font\ninei=cmmi9
    \font\ninesy=cmsy9 
    \skewchar\ninei='177
    \skewchar\ninesy='60

  \font\twelvemrm=cmr10 at 12pt %% new name
  \font\twelvei=cmmi10 at 12pt
  \font\twelvesy=cmsy10 at 12pt
 % \font\twelveex=cmex10 at 12pt

  \font\titlemrm=cmr10 at 14.4pt %% new name
  \font\titlei=cmmi10 at 14.4pt
  \font\titlesy=cmsy10 at 14.4pt
 % \font\titleex=cmex10 at 14.4pt

 %%%% Miscellanious font definitions

  \def\Smallfonts{\ninepoint}

  \def\Hfont{\titlepoint\bf}
  \def\Authorfont{\twelvepoint\it}
  \def\HHfont{\twelvepoint\bf}
  \def\HHHfont{\bf}
  % automatically smaller in 9 point parts
  \def\Bibfont{\tenbf}
  \def\Coordfont{\nineit }% defined in osuPSfnt.sty

  \def \thfont {\bf }
  \def \pffont {\it\itSpacing }
  \def \rkfont {\bf }
  \def \dffont {\bf }
  \def \egfont {\bf }

 %%%%% NINEPOINT %%%%%
 \def\ninepoint{%
  \def\rm{\fam0\ninerm}%
    \textfont0=\ninemrm  \scriptfont0=\sevenrm  \scriptscriptfont0=\fiverm
    \textfont1=\ninei    \scriptfont1=\seveni   \scriptscriptfont1=\fivei
  \def\mit{\fam1\ninei}%
  \def\oldstyle{\fam1\ninei}%
    \textfont2=\ninesy   \scriptfont2=\sevensy  \scriptscriptfont2=\fivesy
    \textfont3=\tenex    \scriptfont3=\tenex    \scriptscriptfont3=\tenex
  \def\it{\fam\itfam\nineit}%
    \textfont\itfam=\nineit
  \def\bf{\ifmmode\fam\bffam\else\ninebf\fi}%
    \textfont\bffam=\ninebold 
    \scriptfont\bffam=\sevenbold 
    \scriptscriptfont\bffam=\fivebold%
  \def\msa{\fam\msafam\ninemsa}%
    \textfont\msafam=\ninemsa 
    \scriptfont\msafam=\sevenmsa
    \scriptscriptfont\msafam=\fivemsa%
  \def\msb{\fam\msbfam\ninemsb}%
    \textfont\msbfam=\ninemsb%
    \scriptfont\msbfam=\sevenmsb%
    \scriptscriptfont\msbfam=\fivemsb%
  \def\eufm{\fam\eufmfam\nineeufm}%
    \textfont\eufmfam=\nineeufm
    \scriptfont\eufmfam=\seveneufm
    \scriptscriptfont\eufmfam=\fiveeufm
   \def\eusm{\fam\eusmfam\nineeusm}%
     \textfont\eusmfam=\nineeusm
     \scriptfont\eusmfam=\seveneusm
     \scriptscriptfont\eusmfam=\fiveeusm
   \def\cyr{\fam\cyrfam\ninecyr}%
     \textfont\cyrfam=\ninecyr
     \scriptfont\cyrfam=\sevencyr
     \scriptscriptfont\cyrfam=\sixcyr%%
  \setbox\strutbox=\hbox{\vrule
      height7pt depth3pt width0pt}%
   \baselineskip=10.8pt\rm}

 \let\eightpoint\ninepoint % we do not use eightpoint

 %%%%% FONTS AT TENPOINT %%%%%
 \def\tenpoint{%
  \def\rm{\fam0\tenrm}%
    \textfont0=\tenmrm \scriptfont0=\sevenrm \scriptscriptfont0=\fiverm%
  \def\mit{\fam1\teni}%
  \def\oldstyle{\fam1\teni}%
    \textfont1=\teni   \scriptfont1=\seveni  \scriptscriptfont1=\fivei%
    \textfont2=\tensy  \scriptfont2=\sevensy \scriptscriptfont2=\fivesy%
    \textfont3=\tenex  \scriptfont3=\tenex   \scriptscriptfont3=\tenex%
  \def\it{\fam\itfam\tenit}%
    \textfont\itfam=\tenit%
  \def\bf{\ifmmode\fam\bffam\else\tenbf\fi}%
    \textfont\bffam=\tenbold% was tenbold for osu
    \scriptfont\bffam=\sevenbold%
    \scriptscriptfont\bffam=\fivebold%
  \def\msa{\fam\msafam\tenmsa}%
    \textfont\msafam=\tenmsa%
    \scriptfont\msafam=\sevenmsa%
    \scriptscriptfont\msafam=\fivemsa%
  \def\msb{\fam\msbfam\tenmsb}%
    \textfont\msbfam=\tenmsb%
    \scriptfont\msbfam=\sevenmsb%
    \scriptscriptfont\msbfam=\fivemsb%
  \def\eufm{\fam\eufmfam\teneufm}%
   \textfont\eufmfam=\teneufm
   \scriptfont\eufmfam=\seveneufm
   \scriptscriptfont\eufmfam=\fiveeufm
   \def\eusm{\fam\eusmfam\teneusm}%
    \textfont\eusmfam=\teneusm
    \scriptfont\eusmfam=\seveneusm
    \scriptscriptfont\eusmfam=\fiveeusm
   \def\cyr{\fam\cyrfam\tencyr}%
    \textfont\cyrfam=\tencyr
    \scriptfont\cyrfam=\sevencyr
    \scriptscriptfont\cyrfam=\sixcyr%%
  \setbox\strutbox=\hbox{\vrule %
      height8.5pt depth3.5ptwidth0pt}%
  \baselineskip=\StdBaselineskip\rm}

 %%%%% FONTS AT TWELVEPOINT %%%%%
 \def\twelvepoint{%
  \def\rm{\fam0\twelverm}%
    \textfont0=\twelvemrm \scriptfont0=\tenmrm \scriptscriptfont0=\sevenrm
    \textfont1=\twelvei   \scriptfont1=\teni   \scriptscriptfont1=\seveni
  \def\mit{\fam1\twelvei}%
  \def\oldstyle{\fam1\twelvei}%
    \textfont2=\twelvesy  \scriptfont2=\tensy  \scriptscriptfont2=\sevensy
    \textfont3=\tenex  \scriptfont3=\tenex  \scriptscriptfont3=\tenex
  \def\it{\fam\itfam\twelveit}%
    \textfont\itfam=\twelveit
  \def\bf{\ifmmode\fam\bffam\else\twelvebf\fi}%
    \textfont\bffam=\twelvebold
    \scriptfont\bffam=\tenbold%
    \scriptscriptfont\bffam=\sevenbold%
  \def\msa{\fam\msafam\twelvemsa}%
    \textfont\msafam=\twelvemsa%
    \scriptfont\msafam=\tenmsa%
    \scriptscriptfont\msafam=\sevenmsa%
  \def\msb{\fam\msbfam\twelvemsb}%
    \textfont\msbfam=\twelvemsb%
    \scriptfont\msbfam=\tenmsb%
    \scriptscriptfont\msbfam=\sevenmsb%
  \def\eufm{\fam\eufmfam\twelveeufm}%
   \textfont\eufmfam=\twelveeufm
   \scriptfont\eufmfam=\teneufm
   \scriptscriptfont\eufmfam=\seveneufm
   \def\eusm{\fam\eusmfam\twelveeusm}%
    \textfont\eusmfam=\twelveeusm
    \scriptfont\eusmfam=\teneusm
    \scriptscriptfont\eusmfam=\seveneusm
   \def\cyr{\fam\cyrfam\tencyr}%
    \textfont\cyrfam=\twelvecyr
    \scriptfont\cyrfam=\tencyr
    \scriptscriptfont\cyrfam=\sevencyr%%
  \setbox\strutbox=\hbox{\vrule
      height10.2pt depth4.55pt width0pt}%
  \baselineskip=14pt\rm}

 %%%%% FONTS AT TITLEPOINT %%%%%
 \def\titlepoint{%
    \textfont0=\titlemrm \scriptfont0=\twelvemrm \scriptscriptfont0=\tenmrm
    \textfont1=\titlei   \scriptfont1=\twelvei   \scriptscriptfont1=\teni
  \def\mit{\fam1\titlei}%
  \def\oldstyle{\fam1\titlei}%
    \textfont2=\titlesy  \scriptfont2=\twelvesy  \scriptscriptfont2=\tensy
    \textfont3=\tenex% math ext not avail in varying sizes??
    \scriptfont3=\tenex
    \scriptscriptfont3=\tenex
  \def\it{\fam\itfam\titleit}%
    \textfont\itfam=\titleit
  \def\bf{\ifmmode\fam\bffam\else\titlebf\fi}%
    \textfont\bffam=\titlebold
    \scriptfont\bffam=\twelvebold%
    \scriptscriptfont\bffam=\tenbold%
  \def\msa{\fam\msafam\titlemsa}%
    \textfont\msafam=\titlemsa%
    \scriptfont\msafam=\twelvemsa%
    \scriptscriptfont\msafam=\tenmsa%
  \def\msb{\fam\msbfam\titlemsb}%
    \textfont\msbfam=\titlemsb%
    \scriptfont\msbfam=\twelvemsb%
    \scriptscriptfont\msbfam=\tenmsb%
  \def\eufm{\fam\eufmfam\titleeufm}%
    \textfont\eufmfam=\titleeufm
    \scriptfont\eufmfam=\twelveeufm
    \scriptscriptfont\eufmfam=\teneufm
   \def\eusm{\fam\eusmfam\titleeusm}%
     \textfont\eusmfam=\titleeusm
     \scriptfont\eusmfam=\twelveeusm
     \scriptscriptfont\eusmfam=\teneusm
   \def\cyr{\fam\cyrfam\tencyr}%
    \textfont\cyrfam=\titlecyr
    \scriptfont\cyrfam=\twelvecyr
    \scriptscriptfont\cyrfam=\tencyr%%
  \setbox\strutbox=\hbox{\vrule
      height12.3pt depth5.54pt width0pt}%
  \baselineskip=16pt\rm}

 %%%% RUNNING HEADINGS
\newbox\AuthorBox\newbox\TitleBox
\newbox\TFLinebox
\newbox\FLinebox
\newbox\HLinebox
\def\SetTFLinebox#1{\setbox\TFLinebox=\hbox{#1}}
\def\SetFLinebox#1{\setbox\FLinebox=\hbox{#1}}
\def\SetHLinebox#1{\setbox\HLinebox=\hbox{#1}}

 \def\SetAuthorHead#1{%
     \setbox\AuthorBox=\hbox{\ninepoint \it 
           \ignorespaces\frenchspacing#1\unskip}}
 \def\SetTitleHead#1{%
     \setbox\TitleBox=\hbox{\ninepoint \it
           \ignorespaces\frenchspacing#1\unskip}}

 %% Italic Spacing Correction
  \def\itSpacing{\relax}
  \def\itSpacingOff{\relax}

  %% Main section headings

 \def\Hrule{\hrule width0pt height0pt}

 %% skip used around proclamations, after section headings,
  % and before subsection-headings:
  \newskip\ProcSkip \ProcSkip 8pt plus2pt minus2pt

 \newskip\LastSkip
 \def\SaveLastSkip{\LastSkip\lastskip}
 \def\RestoreLastSkip{\vskip-\LastSkip\vskip\LastSkip}

 %% Do not indent next paragraph after a header:
 \def\NoindentAfter{\everypar={\setbox0=\lastbox\everypar={}}}

 \long\def\H#1\par#2\par{\notenumber=0 \titlepagetrue%
    {
    \baselineskip=20pt
    \parindent=0pt\parskip=0pt\frenchspacing
    \leftskip=0pt plus .2\hsize minus .3\hsize
    \rightskip=0pt plus .2\hsize minus .3\hsize
 \def\\{\unskip\break}%
    \pretolerance=10000 \Hfont #1\unskip\break
     \vskip7pt\Hrule
\hfill \Authorfont #2\hfill\hfill\unskip}
    \vskip48pt plus 4pt minus 4pt% 60pt=48+12pt
    \par\NoindentAfter\rm}

 \long\def\Hi#1\par#2\par{\notenumber=0 \titlepagetrue%
    {  \baselineskip=0pt  \parindent=0pt\parskip=0pt\frenchspacing
    \leftskip=0pt plus .2\hsize minus .3\hsize
    \rightskip=0pt plus .2\hsize minus .3\hsize
}
    \rm}

 %%% Minor section headings

 \newdimen\PageRemainder
  \def\SetPageRemainder{%\maxdimen case at page tops 12-91 LS
     \PageRemainder=\pagegoal
     \ifdim\PageRemainder=\maxdimen\PageRemainder=\vsize
     \else\advance\PageRemainder by -1\pagetotal\fi}

  \def\Rpt@{}\def\Rpt@@{}

  \long\def\HH#1\par{\par%A
  \SaveLastSkip\removelastskip\goodbreak
  \ifdim\LastSkip<30pt %24pt
     \LastSkip 30pt%24pt 
plus 3pt minus 2pt\fi
  \SetPageRemainder\advance\PageRemainder-\LastSkip
  \ifdim\PageRemainder<150pt
       \edef\Rpt@{remain = \the\PageRemainder\noexpand\\
                pagetotal=\the\pagetotal\noexpand\\
                           pagegoal=\the\pagegoal}%
          \fi
   \ifdim\PageRemainder<65pt %%Head plus 4 lines (approx)
       \ifdim\PageRemainder > 0pt
          \edef\Rpt@@{\noexpand\\
                      Had HH PageRemainder$<$\relax 65pt\noexpand\\
                      Hence forced break!}%
     \vskip 0pt plus .2\PageRemainder\eject %% Pull it out a bit
    \fi\fi
    \vskip\LastSkip\Hrule %%%%%%%%\Hrule added
    \pretolerance=10000\rightskip=0pt plus 3em%B
    \hangafter1 \hangindent=2.2em%
    \noindent
    \HHfont \unskip \Ednote{\Rpt@\Rpt@@}%
            \def\Rpt@{}\def\Rpt@@{}%
            \ignorespaces
            #1\par\rightskip=0pt\pretolerance=\StdPretolerance%
    \NoindentAfter
\tenpoint\rm%
     \medskip \vskip\ProcSkip}%interlineskip adds 2pt to this

  \long\def\HHH#1\par{\par%
  \SaveLastSkip\removelastskip\goodbreak
  \ifdim\LastSkip<\ProcSkip%
     \LastSkip\ProcSkip\fi
  \SetPageRemainder\advance\PageRemainder-\LastSkip
  \ifdim\PageRemainder<150pt
       \edef\Rpt@{remain = \the\PageRemainder\noexpand\\
                pagetotal=\the\pagetotal\noexpand\\
                           pagegoal=\the\pagegoal}%
       \fi
   \ifdim\PageRemainder<48pt  %% 4 lines
        \ifdim\PageRemainder > 0pt
             \edef\Rpt@@{\noexpand\\
                      Had HHH PageRemainder$<$\relax48pt\noexpand\\
                      Hence forced break!}%
       \vskip 0pt plus .2\PageRemainder\eject %% Pull it out a bit
      \fi\fi
   \vskip\LastSkip\par\noindent
   \HHHfont \unskip\Ednote{\Rpt@\Rpt@@}%
  \def\Rpt@{}\def\Rpt@@{}%
  \ignorespaces
   #1\unskip.\quad\rm\ignorespaces
   \ignorepars}

  \long\def\ignorepars#1\par{\def\Test{#1}%
     \ifx\Test\Empty\def\This{\ignorepars}%
        \else\def\This{\Test\par}\fi
           \This}
  \def\Empty{}

 \def\Abstract#1\par{\bgroup\Smallfonts\narrower\HHH #1\par}
 \def\endAbstract{\par\egroup}

 %%%%% Proclamations %%%%%

 \def\ProcBreak{\par%
    \ifdim\lastskip<8pt%
    \removelastskip%
    \penalty-200\vskip\ProcSkip\fi}

 \def\th#1\par{\ProcBreak \noindent
   {\thfont\ignorespaces
    #1\unskip.}\it\itSpacing\kern.4em\ignorepars}%\everymath{\/}

 \def\endth{\ProcBreak\rm\itSpacingOff }%\everymath{}

  %% the theorem statement will be in italic by default

 \def\pf#1\par{\ProcBreak %
    \noindent\pffont#1\unskip.\rm\itSpacingOff{\kern .7em}\ignorepars}

 \def\endpf{\medskip \ProcBreak } %% \qed is alternative

  %% A Box for the QED
  \def\qedbox{\hbox{\vbox{
    \hrule width0.2cm height0.2pt
    \hbox to 0.2cm{\vrule height 0.2cm width 0.2pt
             \hfil\vrule height0.2cm width 0.2pt}
    \hrule width0.2cm height 0.2pt}\kern1pt}}

  %% Typing in \qed makes the qedbox right justified:
  \def\qed{\ifmmode\qedbox
    \else\unskip\ \hglue0mm\hfill\qedbox\ProcBreak\fi}

  \def \rk #1\par{\ProcBreak
     \noindent{\rkfont\ignorespaces #1\unskip.}%
     \rm\kern.6em\ignorepars}

  \def \endrk {\medskip\ProcBreak }

  \def \df #1\par{\ProcBreak
     \noindent{\dffont\unskip\ignorespaces #1\unskip.}%
     \rm\kern.6em\ignorepars}

  \def \enddf {\medskip\ProcBreak }

  \def \eg #1\par{\ProcBreak
     \noindent\egfont\unskip\ignorespaces #1\unskip.
     \rm\kern.6em\ignorepars}

  \def \endeg {\medskip\ProcBreak }

  \newdimen\Overhang

   \def\MaxTag@#1#2#3#4#5{\setbox0=\hbox{#4\ignorespaces#2\unskip}%
     \dimen0=\wd0\advance\dimen0 by#3
     \ifdim\dimen0<#5\relax\dimen0=#5\fi
     \expandafter\edef\csname #1Hang\endcsname{\the\dimen0}}

 \def\MaxItemTag#1{\MaxTag@{Item}{#1}{.4em}{\ItemStyle}{\parindent}}%
 \def\MaxItemItemTag#1{%
        \MaxTag@{ItemItem}{#1}{.4em}{\ItemItemStyle}{\parindent}}
 \def\MaxNrTag#1{\MaxTag@{Nr}{#1}{.5em}{\NrStyle}{\parindent}}
 \def\MaxReferenceTag#1{%
        \MaxTag@{Reference}{[#1]}{.6em}{\ninerm}{\parindent}}
 \def\MaxFootTag#1{\MaxTag@{Foot}{#1}{.4em}{\ninerm}{\z@}}

  %% \SetOverhang@ will prevent for tag-text collision
  \def\SetOverhang@{\Overhang=.8\dimen0%
     \advance\Overhang by \wd0\relax%nec!
     \ifdim\Overhang>\hangindent\relax%nec!
       \advance\Overhang by .25\dimen0%
       \Ednote{Tag is pushing text.}\osumess{Tag is pushing text.}%
     \else\Overhang=\hangindent
     \fi}

   %%% \Item
   \def\Item#1{\par\noindent
      \hangafter1\hangindent=\ItemHang
      \setbox0=\hbox{\ItemStyle\ignorespaces#1\unskip}%
      \dimen0=.4em\SetOverhang@% dimen0 is extra space
      \rlap{\box0}\kern\Overhang\ignorespaces}

   %%% \ItemItem
   \def\ItemItem#1{\par\noindent
      \hangafter1\hangindent=\ItemItemHang
      \setbox0=\hbox{\ItemItemStyle\ignorespaces#1\unskip}%
      \dimen0=.4em\SetOverhang@
      \advance\hangindent by \ItemHang
      \kern\ItemHang\rlap{\box0}%
      \kern\Overhang\ignorespaces}

  %%%% \Nr Items without hanging indentation
  \def\Nr#1{\par\noindent\hangindent=\NrHang %not really a hang
    \setbox0=\hbox{\NrStyle\ignorespaces#1\unskip}%
    \dimen0=.5em\SetOverhang@% dimen0 is extra space
    \rlap{\box0}\kern\Overhang
    \hangindent=\z@\ignorespaces}

  %%%% Roster (not compulsory)
  %%  endRoster has to remember \lastskip (e.g. from a \qed) through \egroup.
   \newskip\Rosterskip\Rosterskip 1pt plus1pt %% modifiable
   \def\Roster{\par\ifdim\lastskip<\Rosterskip\removelastskip\vskip\Rosterskip\fi
    \bgroup}
   \def\endRoster{\par\global\edef\LastSkip@{\the\lastskip}\removelastskip
       \egroup\penalty-50\LastSkip\LastSkip@\relax
       \ifdim\LastSkip<\Rosterskip\LastSkip\Rosterskip\fi
       \vskip\LastSkip}%%changed Feb/5/92 WN

 %%%%% Emphasis %%%%%

 %%%%% Vertical spacing %%%%%

 %%%%% References %%%%%

 \def\cite#1{%\relaxnext@
    \def\nextiii@##1,##2\end@{{\frenchspacing\rm 
      \lBr\ignorespaces##1\unskip{\rm,~\ignorespaces##2}\rBr}}%
    \IN@0,@#1@%
    \ifIN@\def\next{\nextiii@#1\end@}\else
    \def\next{{\rm\lBr#1\rBr}}\fi\next}

 %%%%% Bibliography %%%%%

   \def \Bib#1\par{%
       \par\removelastskip\SetPageRemainder
       \ifdim\PageRemainder < 97pt
        \ifdim\PageRemainder > 0pt
        \vfill\eject
       \fi\fi
    \ProcBreak \par\begingroup\parskip=0 pt%
    \goodbreak \vskip 15 pt plus 10 pt
    \noindent\null\hfill\Bibfont% \kern??pt%  (center over what?)
      \ignorespaces #1\unskip\hfill\null\par 
    \frenchspacing \Smallfonts\rm
    \parskip=2.5 pt plus 1 pt minus.5pt%
    \nobreak\vskip 12pt plus 2pt minus2pt\nobreak
    \leftskip=0 pt \baselineskip=10.5pt}

 \def\ReferenceTagSlide{0em}
  \def\ReferenceTagGap{.5em}

  \def \rf#1{\par\noindent
     \hangafter1\hangindent=\ReferenceHang      
     \setbox0=\hbox{\ninerm[\ignorespaces#1\unskip]}%        
     \dimen0=\ReferenceTagGap\SetOverhang@
     \rlap{\kern\ReferenceTagSlide\box0}%       
     \kern\Overhang\ignorespaces}

  \def\ref#1\par#2\par#3\par#4\par{%
     \rf{#1}#2\unskip,\ #3\unskip,\
     #4\unskip.}

  \def\endBib{\par\endgroup\vskip 12pt minus 6pt }

 %%%%% Coordinates %%%%%

  \long\def\Coordinates#1\endCoordinates{%\relax}
 {\par\vskip4pt\def\\{\unskip, }\Coordfont\baselineskip10.5pt\noindent#1}}

 \def\pagecontents{%\TRMargIns new, \Pagetot@l
  \gdef\Pagetot@l{\pagetotal}
  \ifvoid\TRMargIns\else
    \rlap{\kern\hsize\kern10pt\vbox to 0pt{%
         \box\TRMargIns\vss}}\fi
  \ifvoid\topins\else\unvbox\topins\fi
   \dimen@=\dp\@cclv \unvbox\@cclv % open up \box255
   \ifvoid\footins\else % footnote info is present
     \vskip\skip\footins
     \footnoterule
     \unvbox\footins\fi
   \ifr@ggedbottom \kern-\dimen@ \vfil \fi}

  %%%%% Some math accents %%%%%

 \newcount\Ht %pg121; Height register, used in Linefigure & accents

 \def \Acc{\expandafter } %%% What is this for?? WN

 \def\swthat{\raise -1.1 ex\hbox{\sam$\widehat{}$}}
 \def\swttilde{\raise -1.2 ex\hbox{\sam$\widetilde{}$}}
 \def \overdot{{\raise .2 ex \hbox to 0pt {\hss\bf\smash{.}\hss}}}
 \def \overcircle{{\raise .1 ex \hbox to 0pt
    {\sam$\eightpoint\scriptstyle\hss\circ\hss$}}}

 \def \Mathaccent#1#2{{\sam % E.g. #1=\widehat
  \setbox4=\hbox{$\vphantom{#2}$}
  \Ht=\ht4 %pg120
  \setbox5=\hbox{${#1}$}
  \setbox6=\hbox{${#2}$}
  \setbox7=\hbox to .5\wd6{}
  \copy7\kern .1\Ht \raise\Ht sp\hbox{\copy5}\kern-.1\Ht
  \copy7\llap{\box6}
  }}

  \def\SwtCheck #1{
        \ifmmode \check{#1}%
                \else \v {#1}%
                \fi}

 %%  \barpartial : bar over partial is common, tailor!
 \def\barpartial {%
   \kern .17 em
    \overline {\kern -.17 em\partial\kern-.03 em}%
    \kern .03 em}

 %%%   BEtter overline
 
  \def\Overline#1{\setbox1=\hbox{\sam ${#1}$}%
      \ifdim \wd1 > 6pt
    \kern .11 em
    \overline {\kern -.11 em#1\kern-.14 em}
    \kern .14 em
  \else
    \kern .03 em
    \overline {\kern -.03 em#1\kern-.04 em}
    \kern .04 em
  \fi}

 \def\SOverline#1{\setbox1=\hbox{\sam ${#1}$}%
      \ifdim \wd1 > 7pt
    \kern .22 em
    \overline {\kern -.22 em#1\kern-.09 em}%
    \kern .09 em
  \else
    \kern .10 em
    \overline {\kern -.10 em#1\kern-.04 em}%
    \kern .04 em
  \fi}

  %%% Better underline

 \def\Underline#1{\setbox1=\hbox{\sam ${#1}$}%
      \ifdim \wd1 > 6pt
    \kern .11 em
    \underline {\kern -.11 em#1\kern-.14 em}
    \kern .14 em
  \else
    \kern .03 em
    \underline {\kern -.03 em#1\kern-.04 em}
    \kern .04 em
  \fi}

 \def\SUnderline#1{\setbox1=\hbox{\sam ${#1}$}%
      \ifdim \wd1 > 7pt
    \kern .04 em
    \underline {\kern -.04 em#1\kern-.2 em}%
    \kern .2 em
  \else
    \kern .0 em
    \underline {\kern -.0 em#1\kern-.15 em}%
    \kern .15 em
  \fi}

  %%%%% Miscellaneous %%%%%

 \def \Blackbox
   {\leavevmode\hskip .3pt \vbox
   {\hrule height 5pt\hbox{\hskip 4.5pt}}\hskip .5pt}

 \def \XX{\Blackbox\kern.5pt\Blackbox} %% editorial use

  \def\.{.\kern1pt}

  %% unbreakable hyphen (by local change of hyphenchar to -1)
    \def\Hyphen{\edef\this{\the\hyphenchar\font}%
          \hyphenchar\font=-1\char\this\hyphenchar\font=\this}

  %% Prose In Math or Display 
 \ifx\undefined\text
  \def\text#1{\hbox{\rm #1}}\fi %% AMSTeX is more sophisticated

  %% Math Object Names (multi-character math object names)
  %%\nolimits can be cancelled
                                     % by a following \limits if wanted

%%%% Larry's mathsurround hacks:

   \everymath{}  %% initially, but later ...

  \def\PassMath@@{\aftergroup\AfterMath@} %% use \aftergroup LS 5-92

 \let\PassMath@\PassMath@@

 \def\AfterMath@{\futurelet\next\AfterMathMole@}

 \def\AfterMathMole@{%\show\next
      \ifcat\next\space% picks off CR and \par cases too; not \dots
          \def\this{}%{(space)}%
      \else
      \ifcat\next\egroup %
        \def\this{\osumess{Handset mathsurround?? ---(see dollar brace)}}%
      \else
      \def\this{\AAfterMath@}% this minority case slow
      \fi\fi
      \this}

 \def\hyphen@{-}
 \def\paren@{)}
 \def\apostr@{'}

 \def\MSC#1{\kern-.8\mathsurround#1\kern.8\mathsurround}

 \def\AAfterMath@#1{\def\Next{#1}%\show\Next%
    \IN@0\Next @,.;:!?\relax @%
    \ifIN@\def\this{\MSC{\Next}}%
    \else
    \ifx\Next\hyphen@\def\this{\futurelet\next\AfterHyphen@}%
    \else
    \ifx\Next\paren@\def\this{#1}%
    \else 
    \ifx\Next\apostr@\def\this{#1}%
    \else \def\this{\osumess{Handset mathsurround??}%
                 #1}\fi\fi\fi\fi
    \this}

 \def\AfterHyphen@#1{\def\Next{#1}%
   \ifx\Next\hyphen@\def\this{--}\else
   \ifcat\next\space%
   \def\this{\kern-\mathsurround\kern.05em- \Next}\else
   \def\this{\kern-\mathsurround\kern.05em\Hyphen\Next}\fi\fi\this}

%%%% switches
 \def\sam{\mathsurround=\z@\let\PassMath@\relax}  %
 \def\mas{\mathsurround=\StdMathsurround\let\PassMath@\PassMath@@}
 
 \def\Mas{\mathsurround=\StdMathsurround
                \everymath{\PassMath@}\let\PassMath@\PassMath@@}

 \def\m@th{\mathsurround=\z@\everymath{}}%% good general measure

 \def\m@@th{\mathsurround=\z@\everymath={}\let\m@th\relax}

\def\underbar#1{$\setbox\z@\hbox{#1}\dp\z@\z@
      \m@th \underline{\box\z@}$\relax}

\def\mathhexbox#1#2#3{\leavevmode
  \hbox{\m@@th$\m@th \mathchar"#1#2#3$}}

\def\dots{\relax\ifmmode\ldots\else$\m@th\ldots\,$\relax\fi}
   %%% this first \relax is ONLY original

\def\dotfill{\cleaders\hbox{\m@@th$\m@th \mkern1.5mu.\mkern1.5mu$}\hfill}
\def\rightarrowfill{$\m@th\mathord-\mkern-6mu%
  \cleaders\hbox{\m@@th$\mkern-2mu\mathord-\mkern-2mu$}\hfill
  \mkern-6mu\mathord\rightarrow$\relax}
\def\leftarrowfill{$\m@th\mathord\leftarrow\mkern-6mu%
  \cleaders\hbox{\m@@th$\mkern-2mu\mathord-\mkern-2mu$}\hfill
  \mkern-6mu\mathord-$\relax}

\def\downbracefill{$\m@th\braceld\leaders\vrule\hfill\braceru
  \bracelu\leaders\vrule\hfill\bracerd$\relax}
\def\upbracefill{$\m@th\bracelu\leaders\vrule\hfill\bracerd
  \braceld\leaders\vrule\hfill\braceru$\relax}

\def\angle{{\vbox{\m@@th\ialign{$\m@th\scriptstyle##$\crcr
      \not\mathrel{\mkern14mu}\crcr
      \noalign{\nointerlineskip}
      \mkern2.5mu\leaders\hrule height.34pt\hfill\mkern2.5mu\crcr}}}}

\def\big#1{{\m@@th\hbox{$\left#1\vbox to8.5\p@{}\right.\n@space$}}}
\def\Big#1{{\m@@th\hbox{$\left#1\vbox to11.5\p@{}\right.\n@space$}}}
\def\bigg#1{{\m@@th\hbox{$\left#1\vbox to14.5\p@{}\right.\n@space$}}}
\def\Bigg#1{{\m@@th\hbox{$\left#1\vbox to17.5\p@{}\right.\n@space$}}}
\def\n@space{\nulldelimiterspace\z@ \m@th}

\def\root#1\of{\setbox\rootbox\hbox{\m@@th$\m@th\scriptscriptstyle{#1}$}
  \mathpalette\r@@t}
\def\r@@t#1#2{\setbox\z@\hbox{\m@@th$\m@th#1\sqrt{#2}$\relax}
  \dimen@\ht\z@ \advance\dimen@-\dp\z@
  \mkern5mu\raise.6\dimen@\copy\rootbox \mkern-10mu \box\z@}

\def\mathph@nt#1#2{\setbox\z@\hbox{\m@@th$\m@th#1{#2}$}\finph@nt}

\def\mathsm@sh#1#2{\setbox\z@\hbox{\m@@th$\m@th#1{#2}$}\finsm@sh}

\def\@vereq#1#2{\lower.5\p@\vbox{\m@@th\baselineskip\z@skip\lineskip-.5\p@
    \ialign{$\m@th#1\hfil##\hfil$\crcr#2\crcr=\crcr}}}

\def\mathpalette#1#2{\sam\mathchoice{#1\displaystyle{#2}}%
  {#1\textstyle{#2}}{#1\scriptstyle{#2}}{#1\scriptscriptstyle{#2}}\mas}

\def\widehat#1{\setbox\z@\hbox{\sam$#1$}%
 \ifdim\wd\z@>\tw@ em\mathaccent"0\msbfam@5B{#1}%
 \else\mathaccent"0362{#1}\fi}
\def\widetilde#1{\setbox\z@\hbox{\sam$#1$}%
 \ifdim\wd\z@>\tw@ em\mathaccent"0\msbfam@5D{#1}%
 \else\mathaccent"0365{#1}\fi}

 \def\dots{\relax{}
  \ifmmode\def\thedots{\mdots@}\else\def\thedots{\tdots@}\fi %
  \thedots}

 %% \eqno and \leqno need protection
 \let\@ldeqno\eqno\let\@ldleqno\leqno
 \def\eqno{\everymath{}\@ldeqno} \def\leqno{\everymath{}\@ldleqno}

  \let\@ldeqalignno\eqalignno
  \def\eqalignno#1{\sam\@ldeqalignno{#1}\mas}
  \let\@ldeqalign\eqalign
  \def\eqalign#1{\sam\@ldeqalign{#1}\mas}

 \def\overrightarrow#1{\vbox{\m@th\ialign{##\crcr
      \rightarrowfill\crcr\noalign{\kern-\p@\nointerlineskip}
      $\hfil\displaystyle{#1}\hfil$\crcr}}}
 \def\overleftarrow#1{\vbox{\m@th\ialign{##\crcr
      \leftarrowfill\crcr\noalign{\kern-\p@\nointerlineskip}
      $\hfil\displaystyle{#1}\hfil$\crcr}}}
 \def\overbrace#1{\mathop{\vbox{\m@th\ialign{##\crcr\noalign{\kern3\p@}
      \downbracefill\crcr\noalign{\kern3\p@\nointerlineskip}
      $\hfil\displaystyle{#1}\hfil$\crcr}}}\limits}
 \def\underbrace#1{\mathop{\vtop{\m@th\ialign{##\crcr
      $\hfil\displaystyle{#1}\hfil$\crcr\noalign{\kern3\p@\nointerlineskip}
      \upbracefill\crcr\noalign{\kern3\p@}}}}\limits}

  \let\@ldmatrix\matrix
  \let\end@ldmatrix\endmatrix
  \def\matrix{\sam\@ldmatrix}
  \def\endmatrix{\end@ldmatrix\mas}
  \let\@ldgather\gather
  \let\end@ldgather\endgather
  \def\gather{\sam\@ldgather}
  \def\endgather{\end@ldgather\mas}
  \let\@ldalign\align
  \let\end@ldalign\endalign
  \def\align{\sam\@ldalign}
  \def\endalign{\end@ldalign\mas}
  \let\@ldaligned\aligned
  \let\end@ldaligned\endaligned
  \def\aligned{\sam\@ldaligned}
  \def\endaligned{\end@ldaligned\mas}
  \let\@ldtag\tag
  \def\tag{\sam\@ldtag}
   %
  %%% Commutative diagrams : use LamsCD too?

   \let\MinCDArrowWidth\minCDaw@

  %% will be redefined by BoxedEPS.tex

  %%%%% \FigureTitle %%%%%

%%%% End of Larry's mathsurround stuff
%%%% Start of Walter's insert corrections

\newskip\insertskipamount\newskip\inserthardskipamount
\insertskipamount 6pt plus2pt %This is medskipamount without shrink
\inserthardskipamount 6pt
\def\insertskip{\vskip\insertskipamount}
\newcount\SplitTest%        will be set to -1 if a topinsert has split
\def\SetSplitTest{\SplitTest\insertpenalties
  \insert\topins{\floatingpenalty1}%
  \advance\SplitTest-\insertpenalties}
\def\midinsert{\par
 \SaveLastSkip\penalty-150\SetSplitTest\RestoreLastSkip
 \ifnum\SplitTest=-1
  \@midfalse\p@gefalse\else\@midtrue\fi\@ins}
\def\@ins{\par\begingroup\setbox\z@\vbox\bgroup%
  \vglue\inserthardskipamount}
\def\endinsert{\egroup % finish the \vbox
  \if@mid \dimen@\ht\z@ \advance\dimen@\dp\z@
    \advance\dimen@\insertskipamount%            was 12pt (wn)
    \advance\dimen@\pagetotal\advance\dimen@-\pageshrink
    \ifdim\dimen@>\pagegoal\@midfalse\p@gefalse\fi\fi
  \if@mid%
    \ifdim\lastskip<\insertskipamount\removelastskip\insertskip\fi
    \nointerlineskip\box\z@\penalty-200\insertskip
  \else%
    \SaveLastSkip%                                  added (wn)
    \insert\topins{\penalty100 % floating insertion
    \splittopskip\z@skip
    \splitmaxdepth\maxdimen \floatingpenalty\z@
    \ifp@ge \dimen@\dp\z@
    \vbox to\vsize{\unvbox\z@\kern-\dimen@}% depth is zero
    \else \box\z@\nobreak\insertskip\fi}% was \bigskip\fi (wn)
    \RestoreLastSkip%                               added (wn)
   \fi\endgroup}
%% End Walter's insert stuff

 %%%%% Footnotes %%%%%

  \newcount\notenumber
  
  \def\note{\advance\notenumber by 1
    \footnote{\the\notenumber)}}

  \newbox\footbox

 %% The following modifies Plain TeX definitions, qv
  \def\footnote#1{\let\@sf\empty
    %{(the text)} is read later
    \ifhmode\edef\@sf{\spacefactor\the\spacefactor}\/\fi
    \sam${}^{\fam0 #1}$\@sf\vfootnote{#1}}%

  \def\vfootnote#1{\insert\footins\bgroup
     \interlinepenalty100 \splittopskip=1pt
     \floatingpenalty=20000
     \leftskip=0pt\rightskip=0pt%
     \parindent=.3em%% adjust
     \Smallfonts\rm%%osudeG added \Smallfonts
     \FootItem@{#1}%\strut% not nec
     \futurelet\next\fo@t}

  \def\FootItem@#1{\par\hangafter1\hangindent=\FootHang
     \setbox0=\hbox{\ignorespaces#1\unskip}%
     \dimen0=.4em\SetOverhang@% dimen0 is extra space
     \noindent\rlap{\box0}\kern\Overhang\ignorespaces}

  %\MaxFootTag{2)}%% in param file

  \def\fo@t{\ifcat\bgroup\noexpand\next \let\next\f@@t
    \else\let\next\f@t\fi \next}
  \def\f@@t{\bgroup\aftergroup\@foot\let\next}
  \def\f@t#1{\baselineskip=10pt\lineskip=1pt
            \lineskiplimit=0pt #1\@foot}%
     %%osudeG added \baselineskip=? pt\lineskiplimit=0pt
  \def\@foot{%%% special strut osu for end of each note
        \hbox{\vrule height0pt depth5pt width0pt}
        \egroup}
  \skip\footins=12 pt plus 0pt minus 0pt %% was \bigskipamount
    %% space added when footnote is present
  \count\footins=1000 % footnote magnification factor (1 to 1)
  \dimen\footins=8in % maximum footnotes per page

 %%%% Altenatives

  %%  Editorial stuff (delete??)

 \def\osumess#1{\EdSpider{\immediate\write16{Line \the\inputlineno: #1}}}%
 \def\HideEdStuff{\gdef\EdSpider##1{}}

 \font\BigSym=cmmi10 scaled \magstep 4

 \def\change{\InLMargin{\hbox{\BigSym \char63\kern10pt}}}

 \def\beginchange{\InLMargin{\hbox{\sam\twelvepoint$\heartsuit$\kern10pt}}}

 \def\endchange{\InLMargin{\hbox{\sam\twelvepoint$\spadesuit$\kern10pt}}}

 \def\InLMargin#1{\strut\vadjust{%
     \kern-\strutdepth
     \vtop to \strutdepth{%
         \baselineskip\strutdepth
         \llap{\sam$\smash{\hbox{\EdSpider{#1}}}$}\null}}}

 \def\strutdepth{\dp\strutbox}
 \def\strutheight{\ht\strutbox}

 \def\NoteInRMargin#1{\strut\vadjust{%
     \kern-1.001\strutdepth
     \vtop to \strutdepth{%
       \baselineskip\strutdepth
       \vss\rlap{\ninepoint\unskip\hskip\hsize
         \vtop to 0pt{%
           \hsize=16em\hfuzz=\hsize
           \leftskip=10pt%
           \rightskip=0pt plus 10000pt%
           \baselineskip=9.8pt\lineskip=.2pt%
           \let\\\break
           \noindent\EdSpider{#1}\vss}%
                \kern10pt}\hbox{}}%%\hbox{}=\null crucial!!
       }}

 \def\ednote#1{\NoteInRMargin{\tentt #1}}

 \def\cbar{\InLMargin{%
      \dimen0=\strutdepth\advance\dimen0 by \lineskip
      \vrule width 3pt
      height \strutheight depth \dimen0 \kern
      3pt}}

 \def\ccbar{\InLMargin{%
      \dimen0=2\strutdepth\advance\dimen0 by 2\lineskip
      \vrule width 3pt
        height 3\strutheight depth \dimen0 \kern
      3pt}}

 \newinsert\TRMargIns
 \dimen\TRMargIns=\maxdimen
 %\count\TRMargIns=0
 %\skip\TRMargIns=0pt

  \def\Ednote#1{\insert\TRMargIns{%
       \vbox to 0pt{\hsize=140pt\hfuzz=\hsize
           \leftskip=6pt%
           \rightskip=0pt plus 10000pt%
           \baselineskip=9.8pt\lineskip=.2pt%
           \let\\\break
           %\vglue\pagetotal% misplaces notes if inserts are present
           \SetPageRemainder% This ...
           \vglue540pt\vglue-\PageRemainder%  .. is a fix (WN)
           \noindent\EdSpider{\tentt #1}\vss}%
       \smallskip}}

 \def\KillEdStuff{\def\ednote##1{}\def\Ednote##1{}%
      \let\change\relax\let\beginchange\relax\let\endchange\relax
       \let\cbar\relax\let\ccbar\relax}

 %%% Compatibility with osumrip.sty
  %%

 %%% Parameters
  \topskip=12pt
  \newskip\StdBaselineskip % to set \baselineskip
  \StdBaselineskip 12pt
  \lineskip=1.1pt
  \lineskiplimit=.8pt
  \widowpenalty=10000 % 8000 to 10000
  \clubpenalty=10000  % 8000 to 10000
  \abovedisplayskip=6pt plus 1pt minus 1pt
  \abovedisplayshortskip=3pt plus 1.5pt
  \belowdisplayskip=6pt plus 1pt minus 1pt
  \belowdisplayshortskip=5pt plus 1pt minus 1pt
  \hfuzz=1.5pt   % Enable overfull box warnings at console

  \def\StdPretolerance{100}
  \tolerance=\StdPretolerance

  \newdimen\StdMathsurround
  \StdMathsurround=1.5pt % 1pt usual without \Mas
  \mathsurround=\StdMathsurround
  \Mas                   %% sophisticated mathsurround on
 % \Sam                   %% sophisticated mathsurround off

%% marker before English punctuation in displayed math
   \def\prose{\relax\hbox{\kern.6\StdMathsurround}}
  
  \def\StdParskip{0pt}    %% Larry wants {2pt plus 1pt}
  \parskip=\StdParskip
  \parindent=0.5cm
 
%%%% load Times for main body font

  \def\Times{ptmr  } 
  \def\TimesI{ptmri  } 
  \def\TimesB{ptmb  }
  \def\TimesBI{ptmbi  }
  \def\HelveticaN{phvrrn }

  =\Times at 10bp% roman text
  =\TimesB at 10bp% boldface extended
   % slanted roman
  \font\tenit=\TimesI at 10bp% text italic
  =\TimesBI at 10bp

  \font\tenmrm=cmr10  %%new name for math role at full size

%%%%% Fonts at ninepoint %%%%%

    =\Times at 9bp 
    \font\nineit=\TimesI at 9bp 
    =\TimesB at 9bp 
    =\TimesBI at 9bp 

    =\HelveticaN at 9bp 
       % see below

%%%%% Fonts at twelvepoint %%%%%

  =\Times at 12bp
  \font\twelveit=\TimesI at 12bp
  =\TimesB at 12bp

%%%%% Fonts at titlepoint %%%%%

  \font\titleit=\TimesI at 14.4bp
  =\TimesB at 14.4bp

 \SetAuthorHead{AuthorHead} % needs \ninepoint since box set
 \SetTitleHead{TitleHead}  % notably \HeaderFont

%%%% Char adjustments %%%%

  \def\lBr{\raise.125ex\hbox{[\kern.1125ex}}
  \def\rBr{\raise.125ex\hbox{\kern.1125ex]}}

 \setbox\footbox=\hbox{\Smallfonts 2)~}

%% Some optional font dimension and spacing 
%% adjustments beyond this point

%% Correct the lousy spacing of italic f (a hack).

  \bgroup
  \catcode`\@=11 %localised
  \gdef\itSpacing{%
     \xspaceskip=.31em plus.1em minus.05em \sfcode `f=2001
     \itWarning@\let\itWarning@\itWarning@@}
  \gdef\itSpacingOff{%
     \xspaceskip=0pt \sfcode `f=1000
     \let\itWarning@\relax}
   \global\let\itWarning@\relax
  \gdef\itWarning@@{\errmessage{%
  Special italic spacing already in force
  (you have probably omitted an ``endth'').
  See itSpacing macro in osuPSfnt.sty
         }}
  \egroup

 %%% Provisional fontdimen settings
  %%
 \fontdimen1\titlebf=0.0pt
 \fontdimen2\titlebf=3.6135pt
 \fontdimen3\titlebf=2.8908pt
 \fontdimen4\titlebf=1.44539pt
 \fontdimen5\titlebf=6.64882pt
 \fontdimen6\titlebf=14.45398pt
 \fontdimen7\titlebf=1.60439pt

 \fontdimen1\tenbi=0.26794pt
 \fontdimen2\tenbi=2.50937pt
 \fontdimen3\tenbi=2.00749pt
 \fontdimen4\tenbi=1.00374pt
 \fontdimen5\tenbi=4.59717pt
 \fontdimen6\tenbi=10.03749pt
 \fontdimen7\tenbi=1.11415pt

 \fontdimen1\twelverm=0.0pt
 \fontdimen2\twelverm=3.01125pt
 \fontdimen3\twelverm=2.409pt
 \fontdimen4\twelverm=1.2045pt
 \fontdimen5\twelverm=5.39615pt
 \fontdimen6\twelverm=12.045pt
 \fontdimen7\twelverm=1.33699pt

 \fontdimen1\twelveit=0.27731pt
 \fontdimen2\twelveit=3.01125pt
 \fontdimen3\twelveit=2.409pt
 \fontdimen4\twelveit=1.2045pt
 \fontdimen5\twelveit=5.37207pt
 \fontdimen6\twelveit=12.045pt
 \fontdimen7\twelveit=1.33699pt

 \fontdimen1\twelvebf=0.0pt
 \fontdimen2\twelvebf=3.01125pt
 \fontdimen3\twelvebf=2.409pt
 \fontdimen4\twelvebf=1.2045pt
 \fontdimen5\twelvebf=5.5407pt
 \fontdimen6\twelvebf=12.045pt
 \fontdimen7\twelvebf=1.33699pt

 \fontdimen1\tenrm=0.0pt
 \fontdimen2\tenrm=2.50937pt
 \fontdimen3\tenrm=2.00749pt
 \fontdimen4\tenrm=1.00374pt
 \fontdimen5\tenrm=4.49678pt
 \fontdimen6\tenrm=10.03749pt
 \fontdimen7\tenrm=1.11415pt

 \fontdimen1\tenit=0.27731pt
 \fontdimen2\tenit=2.50937pt
 \fontdimen3\tenit=2.00749pt
 \fontdimen4\tenit=1.00374pt
 \fontdimen5\tenit=4.47672pt
 \fontdimen6\tenit=10.03749pt
 \fontdimen7\tenit=1.11415pt

 \fontdimen1\tenbf=0.0pt
 \fontdimen2\tenbf=2.50937pt
 \fontdimen3\tenbf=2.00749pt
 \fontdimen4\tenbf=1.00374pt
 \fontdimen5\tenbf=4.61723pt
 \fontdimen6\tenbf=10.03749pt
 \fontdimen7\tenbf=1.11415pt

 \fontdimen1\ninerm=0.0pt
 \fontdimen2\ninerm=2.25842pt
 \fontdimen3\ninerm=1.80673pt
 \fontdimen4\ninerm=0.90337pt
 \fontdimen5\ninerm=4.0471pt
 \fontdimen6\ninerm=9.03374pt
 \fontdimen7\ninerm=1.00273pt

 \fontdimen1\nineit=0.27731pt
 \fontdimen2\nineit=2.25842pt
 \fontdimen3\nineit=1.80673pt
 \fontdimen4\nineit=0.90337pt
 \fontdimen5\nineit=4.02904pt
 \fontdimen6\nineit=9.03374pt
 \fontdimen7\nineit=1.00273pt

 \fontdimen1\ninebf=0.0pt
 \fontdimen2\ninebf=2.25842pt
 \fontdimen3\ninebf=1.80673pt
 \fontdimen4\ninebf=0.90337pt
 \fontdimen5\ninebf=4.15552pt
 \fontdimen6\ninebf=9.03374pt
 \fontdimen7\ninebf=1.00273pt

 %%% \SetExtraSpaces \MaxSpaceFactor \SetSpaceFactors
  %%  See TeXbook, page 76.

 \newcount\MaxSpaceFactor
 \MaxSpaceFactor=3000 %% to reset later

 %%%%% Tag styles and (hang-) indents
 \def\ItemStyle{\rm}
 \def\NrStyle{\rm}
 \def\ItemItemStyle{\rm}

 %% Analog dimensioning, convenient for local modifications:
 \MaxItemTag{(iii)}
 \MaxItemItemTag{(iii)}
 \MaxNrTag{(2)}
 \MaxFootTag{2)}
 % \MaxReferenceTag{AaaAA} % for biblio
 \def\ReferenceHang{30pt}

 \catcode`\@=\active

%%%%% End of hack of Neumann-Siebenmann macros

\loadbold

=\Times  
=\Times scaled750
=\Times scaled650
\font\rms=\Times scaled 920 

=\TimesBI scaled 860
=\TimesI scaled 860

\textfont0=\rrm  
\scriptfont0=\erm 
\scriptscriptfont0=\srm

\def\Augment#1#2{%
    \toks0\expandafter{#1}\toks2{#2}%
    \edef#1{\the\toks0\the\toks2}}

 \font\twelverma=\Times  scaled 1200
 \font\tenrma=\Times  scaled 1000
 \font\ninerma=\Times scaled 920
 =\Times scaled 840
 \font\sevenrma=\Times scaled 760
 =\Times scaled 680
 \font\fiverma=\Times scaled 600

 \Augment\tenpoint{%
  \textfont0=\tenrma  \scriptfont0=\sevenrma  
  \scriptscriptfont0=\fiverma  }

 \Augment\ninepoint{%
  \textfont0=\ninerma  \scriptfont0=\sevenrma 
  \scriptscriptfont0=\fiverma}

 \Augment\twelvepoint{%
  \textfont0=\twelverma  \scriptfont0=\ninerma  
  \scriptscriptfont0=\sevenrma}

\mathsurround=1pt
\hsize=13.45truecm
\vsize=19.5truecm
\hoffset=1.25truecm
\voffset=2truecm
\advance\baselineskip by 2pt

\predefine\til{\~}
\def\~#1{\relax\ifmmode\widetilde{#1}\else\til{#1}\fi}

\redefine \le{\leqslant}
\redefine \ge{\geqslant}
\define \wt#1{\mathaccent"0365{#1}}
\define \wh#1{\mathaccent"0362{#1}}

\define \iss{\,\Mathaccent{\raise -.8 ex\hbox{$\widetilde{}$\kern.1em}}\rightarrow\,}

\define \mo{\mathop{\fam0 mod\,\,}}

\define \gcdd{\mathop{\fam0 gcd}}

\define \id{\operatorname{\fam0 id\,}}

\define \chr{\mathop{\fam0 char}\,}

\define \res{\operatorname{\fam0 res}}

\define \Br{\operatorname{\fam0 Br}}

\define \Aut{\operatorname{\fam0 Aut}}

\define \Int{\mathop{\fam0 Int}}

\Mas
\HideEdStuff
\rm 
 
%%%% For GT headers and footers:

\def\issn{{\nineit ISSN 1464-8997 (on line) 1464-8989 (printed)}}

\def\gtp{{\nineit Published 10 December 2000: \ \copyright\ Geometry \& 
Topology Publications}}

\def\gtv3{{\nineit Geometry \& Topology Monographs, Volume 3 (2000) --
Invitation to higher local fields}}

%%%%% For section idents:

\def\lione
{{\rms Geometry \& Topology Monographs}}

\def \litwo{{\rms Volume 3: Invitation to higher local fields
}} 

\def\tinfo #1.#2.#3-#4
{{
\noindent  {\lione} \hfill 
\par 
\vskip-1.5pt
\noindent {\litwo} \hfill
\par 
\vskip-1,5pt
\noindent {\rms Part #1, section #2, pages #3--#4} \hfill
\vskip24pt 
}}

\def\tinfos #1.#2.#3-#4
{{
\noindent  {\lione} \hfill 
\par 
\vskip-1.5pt
\noindent {\litwo} \hfill
\par 
\vskip-1.5pt
\noindent {\rms Pages #3--#4} \hfill
\vskip24pt 
}}

\def\tinfoi #1
{{
\noindent  {\lione} \hfill 
\par 
\vskip-1.5pt
\noindent {\litwo} \hfill
\par 
\vskip-1.5pt
\noindent {\rms Pages iii--xi: Introduction and contents} \hfill
\vskip26pt 
}}

%%%% Set headers and footers %%%%

  \def\titlepagehead{\hfil}

  \newif\iftitlepage\titlepagefalse
  \newif\ifblankpage\blankpagefalse
  \def\makeheadline{
     \ifblankpage{}\else%
     \iftitlepage
\vbox{\line{\vbox to 8.5pt{}
\ninerm
\copy\HLinebox \hfill
\hglue5mm\ninebf\folio 
\titlepagehead}}%
      \else
\vbox{\ifodd\pageno\rightheadline\else\leftheadline\fi}%
      \fi\vskip 12pt\fi}%
     \def\rightheadline{\line{\vbox to 8.5pt{}%
      \ninerm
\copy\TitleBox \hfill
\hglue5mm\ninebf\folio}}%
     \def\leftheadline{\line{\vbox to 8.5pt{}%
        \unskip\ninerm\unskip\ninebf\folio\hglue5mm
      %*%
 \hfill \copy\AuthorBox
%\hfill
}}

 \footline={\ifblankpage{}\else
\iftitlepage\ninepoint\sam\hfill%} 
\line{\vbox to 8.5pt{}%\ninerm
\copy\TFLinebox
\hfill
\hglue5mm %\ninebf\folio
}
            \else
\ninepoint\sam\hfill%}
\line{\vbox to 8.5pt{}%\ninerm
\copy\FLinebox
\hfill 
\hglue5mm
}
\hfil\fi\global\titlepagefalse\fi}

\def\blankpage{{\blankpagetrue\noindent\hbox to 10pt{\hss}\vfill
\pagebreak}}

\tenpoint\rm %% always start here
 
  %%% all done and macros loaded!

\pageno=281

\tinfo II.8.281-292

\SetTFLinebox{\gtp }
\SetFLinebox{\gtv3 }
\SetHLinebox{\issn}

\H 8. Higher local skew fields

Alexander Zheglov

\SetAuthorHead{A. Zheglov}
\SetTitleHead{Part II. Section 8. Higher local skew fields
\qquad\qquad}

$n$-dimensional local skew fields  are a natural generalization of  $n$-dimensional local fields. The latter  have numerous applications to problems of algebraic geometry, both arithmetical and geometrical, as it is shown in this volume. From this viewpoint, it would be reasonable to restrict oneself to commutative fields only. Nevertheless, 
already in  class field theory one meets non-commutative rings
which are  skew fields  finite-dimensional over their center $K$. 
For example,  $K$  is a (commutative) local field and the skew field  represents   elements of the Brauer group of the field $K$ (see also an example below). In \cite{Pa} A.N. Parshin  pointed out  another class of non-commutative local fields arising in differential equations and showed that these skew fields possess many features of  commutative fields. He defined a skew field of formal pseudo-differential operators in $n$ variables and studied some of their properties. He raised a problem of  classifying  non-commutative local skew fields. 

In this section we treat the case of  $n=2$ and list a number of results, in particular a classification of certain types of 2-dimensional local skew fields.

\HH 8.1. Basic definitions

\df Definition

A skew field $K$ is called a {\it complete discrete valuation skew field}
if $K$ is complete with respect to a discrete valuation
(the residue skew field is not necessarily commutative).
A field $K$ is called an {\it $n$-dimensional local skew field}
if there are skew fields
$K=K_n,K_{n-1},\dots,K_0$ such that
each $K_i$ for $i>0$ is a complete discrete valuation skew field
with residue skew field $K_{i-1}$.
\enddf

\eg Examples 

\Roster
\Item{(1)} Let $k$ be a field. Formal pseudo-differential operators
over $k((X))$ form a 2-dimensional local skew field
$K=k((X))((\partial_X^{-1}))$,
$\partial_X X=X\partial_X+1$. If $\chr( k)=0$ we get an example of a skew field which is an infinite dimensional vector space over its centre.

\Item{(2)} Let $L$ be a local field of equal characteristic (of any dimension).
Then an element of $\Br(L)$ is an example of 
a skew field which is finite dimensional over its centre.
\endRoster
\endeg

\smallskip

From now on {\it let} $K$ be a two-dimensional local skew field.
Let $t_2$ be a generator of $\Cal M_{K_2}$
and $t_1'$ be a generator of $\Cal M_{K_1}$.
If  $t_1\in K$ is a lifting of $t_1'$ then 
$t_1,t_2$ is  called a {\it system of local parameters} of $K$. 
We denote by $v_{K_2}$ and $v_{K_1}$ the (surjective) discrete valuations 
of $K_2$ and $K_1$ associated with $t_2$ and $t_1'$.

\df Definition

 A two-dimensional local skew field $K$ is said to {\it split}
if there is a section of the homomorphism
$\Cal O_{K_2}\to K_1$ where $\Cal O_{K_2}$ is the ring of integers of $K_2$.
\enddf

\eg Example (N. Dubrovin)

Let $\Bbb Q\,((u))\langle x,y\rangle$ be a free associative algebra over $\Bbb
Q\,((u))$
with generators $x,y$.
Let $I=\langle [x,[x,y]],[y,[x,y]]\rangle$.
Then the quotient $$A=\Bbb Q\,((u))\langle x,y\rangle/I$$ is a $\Bbb Q$-algebra which
has no non-trivial 
zero divisors, and
 in which $z=[x,y]+I$ is a central element.
Any element of $A$ can be uniquely represented in the form 
$$
f_0+f_1z+\ldots +f_mz^m
$$ 
where $f_0,\ldots ,f_m$ are polynomials in the variables $x,y$.
% and  monomials of the polynomial $f_j$ are ordered:
%$$
%f_j=a+bx+cy+d_1x^2+d_2xy+d_3y^2+\ldots 
%$$
%$a,b,c,d_i\in \Bbb Q\, ((u))$. 

One can define a discrete valuation $w$ on $A$ such that $w(x)=w(y)=w(\Bbb Q\,((u)))=0$, $w([x,y])=1$, $w(a)=k$ if $a=f_kz^k+\ldots +f_mz^m$, $f_k\not=0$. 
The skew field $B$ of fractions of $A$ has a discrete valuation
$v$ which is a unique extension of $w$. 
The completion of $B$ with respect to $v$ is
a two-dimensional local skew field which does not split (for details see \cite{Zh, Lemma 9}).
\endeg

\df Definition

Assume  that $K_1$ is a field.
The homomorphism $$\varphi_0\colon K^*\to \Int(K), \quad \varphi_0(x)(y)=x^{-1}yx$$
induces a homomorphism $\varphi\colon K_2^*/\Cal O_{K_2}^*\to \Aut(K_1)$.
The {\it canonical automorphism} of $K_1$ is $\alpha=\varphi(t_2)$
where $t_2$ is an arbitrary prime element of $K_2$.
\enddf

\df Definition

Two two-dimensional local skew fields $K$ and $K'$ are {\it isomorphic} if
there is an isomorphism $K\to K'$ 
which maps $\Cal O_K$ onto $\Cal O_{K'}$,
$\Cal M_K$ onto $\Cal M_{K'}$
and $\Cal O_{K_1}$ onto $\Cal O_{K_1'}$,
$\Cal M_{K_1}$ onto $\Cal M_{K_1'}$.   
\enddf

\HH 8.2. Canonical automorphisms of infinite order

\th Theorem

\Roster 
\Item{(1)}
Let $K$ be a two-dimensional local skew field. 
 If $\alpha^n\not=\id$ for all $n\ge 1$ then
$\chr(K_2)=\chr(K_1)$,
$K$ splits and $K$ is isomorphic to a two-dimensional local skew field
 $K_1((t_2))$ where $t_2a=\alpha(a)t_2$ for all $a\in K_1$.

\Item{(2)} Let $K,K'$ be two-dimensional local skew fields
and let $K_1, K_1'$ be fields. Let   $\alpha^n\not=\id$, ${\alpha'}^n\not=\id$
for all $n\ge 1$.
Then $K$ is isomorphic to $K'$
if and only if 
there is an isomorphism $f\colon K_1\to K_1'$ such that
$\alpha=f^{-1}\alpha' f$ where $\alpha,\alpha'$ are the canonical
automorphisms of $K_1$ and $K'_1$.
\endRoster 
\endth

\rk Remarks

\Roster
\Item{1.} This theorem is true for any higher local skew field. 

\Item{2.} There are examples (similar to Dubrovin's example)
of local skew fields which do not split and in which 
$\alpha^n= \id$ for some positive integer $n$. 
\endRoster 
\endrk

\pf Proof

(2) follows from (1). We sketch the proof of
(1). For details see  \cite{Zh, Th.1}.

If $\chr (K)\not=\chr
(K_1)$ then $\chr (K_1)=p>0$. Hence $v (p)= r>0$. Then for any element $t\in
K$ with $v (t)=0$ we have $ptp^{-1}\equiv \alpha^r(\Overline{t})\mo \Cal M_K$ where
$\Overline{t}$ is the image of $t$ in $K_1$. But on the other hand, $pt=tp$,
 a contradiction.

Let $F$ be the prime field in $K$. Since $\chr (K)=\chr
(K_1)$ the field $F$ is a subring of $\Cal O=\Cal O_{K_2}$.  
One can easily show  
that there exists an element $c\in K_1$ such that ${\alpha^n}(c)\not= c$ for every  $n\ge 1$ \cite{Zh, Lemma 5}. 

Then any lifting $c'$ in $\Cal O$ of $c$ 
is transcendental over $F$.
 Hence we can embed the field $F
(c')$ in $\Cal O$. Let $\Overline L$ be a maximal field extension of $F (c')$ 
which can be embedded in $\Cal O$. Denote by $L$ its image in $\Cal O$. Take
$\Overline{a}\in K_1\setminus \Overline{L}$. 
We claim that  there exists
a lifting  $a'\in \Cal O$ of $\Overline{a}$ such that $a'$ commutes with
every element in $L$. To prove this fact we use the completeness of $\Cal O$
in the following argument.

Take any lifting $a$ in $\Cal O$ of $\Overline{a}$. For every element $x\in L$  
 we have $axa^{-1}\equiv x\mo \Cal M_K$. If $t_2$  is a prime element of $K_2$
we can write
$$
axa^{-1}=x+\delta_1(x)t_2
$$
where ${\delta_1}(x)\in \Cal O$. 
The map $\Overline{\delta_1}\colon L\ni x\to 
\Overline{\delta_1(x)}\in K_1$ is an $\alpha$-derivation, i.e.
$$\Overline{\delta_1}(ef)=\Overline{\delta_1}(e)\alpha(f)+e\Overline{\delta_1}(f)$$
for all $e,f\in L$. 
Take  an element $h$ such that $\alpha (h)\not=h$,
then $\Overline{\delta_1}(a)=g\alpha(a)-ag$ where
$g=\Overline{\delta_1}(h)/(\alpha(h)-h)$. 
Therefore there is $a_1\in K_1$ such that $$(1+a_1t_2)axa^{-1}(1+a_1t_2)^{-1}\equiv 
x\mo \Cal M_K^2.$$ By induction we can find an element $a'=\dots \cdot (1+a_1t_2)a$ such
that $a'x{a'}^{-1}=x$.

Now, if $\Overline{a}$ is not algebraic over $\Overline{L}$,
then for its lifting $a'\in \Cal O$ which commutes with $L$
we would deduce that $L(a')$ is a  
field extension of $F (c')$ 
which can be embedded in $\Cal O$, which contradicts the maximality of $L$. 

Hence $\Overline{a}$ is algebraic and separable over $\Overline{L}$.
Using a generalization of Hensel's Lemma \cite{Zh, Prop.4} 
we can find a lifting $a'$ of $\Overline{a}$
such that $a'$ commutes with elements of $L$ and
$a'$ is algebraic over $L$, which again leads to a contradiction.

Finally  let $\Overline{a}$ be purely inseparable over $\Overline{L}$,
$\Overline{a}^{p^k}=\Overline{x}$, $x\in L$. 
Let $a'$ be its lifting which commutes with every element of $L$.
Then ${a'}^{p^k}-x$ commutes with every element of $L$.
If $v_K({a'}^{p^k}-x)=r\not=\infty$ then similarly to the beginning of this proof we deduce that the image of 
$({a'}^{p^k}-x)c({a'}^{p^k}-x)^{-1}$ in $K_1$
is equal to $\alpha^r(c)$ (which is distinct from $c$),
a contradiction.
Therefore, ${a'}^{p^k}=x$ and the field 
$L(a')$ is a  
field extension of $F (c')$ 
which can be embedded in $\Cal O$, which contradicts the maximality of~$L$. 

Thus, $\Overline{L}=K_1$.

To prove that $K$ is isomorphic to a skew field $K_1((t_2))$ where $t_2a=\alpha (a)t_2$ one can apply similar arguments as in the proof of the existence of an element $a'$ such that $a'x{a'}^{-1}=x$ (see above). So, one can find a parameter $t_2$ with a given property. 
\qed\endpf

In some cases we have a complete classification of local skew fields.

\th Proposition {{\rm(\cite{Zh})}}

Assume that $K_1$ is isomorphic to $k((t_1))$. Put
$$\zeta =\alpha (t_1)t_1^{-1}\mo \Cal M_{K_1}.$$
Put $i_{\alpha}=1$ if $\zeta$ is not a root of unity in $k$ and
$i_{\alpha}=v_{K_1}(\alpha^n(t_1)-t_1)$ if $\zeta$ is a primitive $n$th root. 
Assume that $k$ is of characteristic zero.  
Then
there is an automorphism $f\in \Aut_{k}(K_1)$ such that
$f^{-1} \alpha f =\beta $
where
$$\beta (t_1)= \zeta t_1+xt_1^{i_{\alpha}}+x^2yt_1^{2i_{\alpha}-1}$$ for some 
$x\in k^*/k^{*(i_{\alpha}-1)}$, $y\in k$.

Two automorphisms $\alpha$ and $\beta$ are conjugate if and only if
 $$(\zeta(\alpha) , i_{\alpha}, x(\alpha), y(\alpha) )=
(\zeta(\beta) , i_{\beta}, x(\beta), y(\beta) ).$$ 
\endth

\pf Proof

First we prove that $\alpha=f\beta' f^{-1}$ where 
$$
\beta' (t_1)=\zeta t_1+xt_1^{in+1}+yt_1^{2in+1}
$$
for some natural $i$. Then we prove that $i_{\alpha}=i_{\beta'}$. 

Consider a set $\{\alpha_i : i\in \Bbb N\,\}$ where $\alpha_i=f_i\alpha_{i-1}f_i^{-1}$, $f_i(t_1)=t_1+x_it_1^i$ for some $x_i\in k$, $\alpha_1=\alpha$. Write 
$$
\alpha_i(t_1)=\zeta t_1+a_{2,i}t_1^2+a_{3,i}t_1^3+\ldots.
$$
One can check that  $a_{2,2}=x_2(\zeta^2-\zeta )+a_{2,1}$ and hence there exists an element $x_2\in k$ such that $a_{2,2}=0$. Since $a_{j,i+1}=a_{j,i}$, we have $a_{2,j}=0$ for all $j\ge 2$. 
Further, $a_{3,3}=x_3(\zeta^3-\zeta )+a_{3,2}$ and hence there exists an element $x_3\in k$ such that $a_{3,3}=0$. Then $a_{3,j}=0$ for all $j\ge 3$. Thus, any element $a_{k,k}$ can be made equal to zero if $n {\not|} (k-1)$, 
and therefore $\alpha =f\tilde{\alpha}f^{-1}$ where 
$$
\tilde{\alpha}(t_1)=\zeta t_1+{\tilde{a}}_{in+1}t_1^{in+1}+{\tilde{a}}_{in+n+1}t_1^{in+n+1}+\ldots
$$
for some $i$, ${\tilde{a}}_j\in k$. Notice that ${\tilde{a}}_{in+1}$ does not depend on $x_i$. Put $x=x(\alpha )={\tilde{a}}_{in+1}$. 

Now we replace $\alpha$ by $\tilde{\alpha}$.
One can check that if $n|(k-1)$ then  
$$
a_{j,k}=a_{j,k-1}\qquad \text{\rm for $2\le j< k+in$} 
$$ 
and 
$$
a_{k+in,k}=x_kx(k-in-1)+a_{k+in}+\text{\rm \ some polynomial which does not depend on } x_k. 
$$
From this fact it immediately follows that $a_{2in+1,in+1} $ does not depend on $x_i$ and for all $k\ne in+1$ $a_{k+in,k}$ can be made equal to zero. 
Then $y=y(\alpha )=a_{2in+1,in+1}$. 

Now we prove that $i_{\alpha}=i_{\beta'}$. Using the formula
$$
{\beta'}^n(t_1)=t_1+nx(\alpha )\zeta^{-1}t_1^{in+1}+\ldots 
$$
we get $i_{\beta'}=in+1$. Then one can check that $v_{K_1}(f^{-1}(\alpha^n-\id)f)=v_{K_1}(\alpha^n-\id)=i_{\alpha}$. Since ${\beta'}^n-\id=f^{-1}(\alpha^n-\id)f$, we get the identity $i_{\alpha}=i_{\beta'}$.

The rest of the proof is clear. For details see \cite{Zh, Lemma 6 and  Prop.5}. \qed\endpf

\HH 8.3. Canonical automorphisms of finite order

\HHH 8.3.1. Characteristic zero case

\phantom{}\par 

Assume that 
\Roster
\Item{} a two-dimensional local skew field $K$ splits, 
\Item{} $K_1$ is a field, $K_0\subset Z(K)$, 
\Item{} $\chr(K)=\chr(K_0)=0$,  
\Item{} $\alpha^n=\id$ for some $n\ge 1$, 
\Item{} for any convergent sequence $(a_j)$ in $K_1$ the sequence $(t_2a_jt_2^{-1})$ converges in $K$.
\endRoster  

\th Lemma

$K$ is isomorphic to a two-dimensional local skew field
$K_1((t_2))$ where
$$
t_2at_2^{-1}=\alpha(a)+{\delta_i}(a)t_2^i+{\delta_{2i}}(a)t_2^{2i}+{\delta_{2i+n}}(a)t_2^{2i+
n}+\ldots \qquad \text{for all  $a\in K_1$}
$$
where $n|i$ and $\delta_j: K_1\to K_1$ are linear maps and
$$
\delta_i(ab) = \delta_i(a)\alpha (b)+\alpha (a)\delta_i(b)\qquad
\text{for every $a,b\in K_1$.} 
$$
Moreover
$$
t_2^nat_2^{-n}=a+\delta'_i(a)t_2^i+\delta'_{2i}(a)t_2^{2i}+\delta'_{2i+n}(a)t_2^{2i+n}+\ldots 
$$
where $\delta'_j$ are linear maps and $\delta'_i$ and $\delta:= \delta'_{2i}-((i+1)/2){\delta'_i}^2$ are derivations.
\endth

\rk Remark

The following fact holds for the field $K$ of any characteristic:
$K$ is isomorphic to a two-dimensional local skew field
$K_1((t_2))$ where
$$
t_2at_2^{-1}=\alpha (a)+\delta_i(a)t_2^i+\delta_{i+1}(a)t_2^{i+1}+\ldots
$$
where $\delta_j$ are linear maps which satisfy some identity. For explicit formulas see \cite{Zh, Prop.2 and Cor.1}.

\pf Proof

It is clear that $K$ is isomorphic to a two-dimensional local skew field $K_1((t_2))$ where 
$$
t_2at_2^{-1}=\alpha (a) +\delta_1(a)t_2+\delta_2(a)t_2^2+\ldots\qquad
\text{\rm for all $a$ } 
$$
and $\delta_j$ are linear maps. Then $\delta_1$ is a $(\alpha^2, \alpha )$-derivation, that is $\delta_1(ab)=\delta_1(a)\alpha^2 (b)+\alpha (a)\delta_1(b)$.

Indeed,
$$
\aligned &t_2abt_2^{-1}=t_2at_2^{-1}t_2bt_2^{-1}=(\alpha (a)+\delta_1(a)t_2+\ldots )(\alpha (b)+\delta_1(b)t_2+\ldots )
\\ 
=&\,\alpha (a)\alpha (b)+(\delta_1(a)\alpha^2(b)+\alpha (a)\delta_1(b))t_2+\ldots=
\alpha (ab)+\delta_1(ab)t_2+\ldots  .
\endaligned 
$$
From the proof of Theorem 8.2 it follows that $\delta_1$ is an inner derivation, i.e. $\delta_1(a)=g\alpha^2(a)-\alpha (a)g$ for some $g\in K_1$, and that there exists  a $t_{2,2}=(1+x_1t_2)t_2$ such that 
$$
t_{2,2}at_{2,2}^{-1}=\alpha (a)+\delta_{2,2}(a)t_{2,2}^2+\ldots .
$$  
One can easily check that $\delta_{2,2}$ is a $(\alpha^3,\alpha )$-derivation. Then it is an inner derivation and there exists $t_{2,3}$ such that 
$$
t_{2,3}at_{2,3}^{-1}=\alpha (a)+\delta_{3,3}(a)t_{2,3}^3+\ldots  .
$$
By induction one deduces that if 
$$
t_{2,j}at_{2,j}^{-1}=\alpha (a)+\delta_{n,j}(a)t_{2,j}^n+\ldots +\delta_{kn,j}(a)t_{2,j}^{kn}+\delta_{j,j}(a)t_{2,j}^j+\ldots
$$
then $\delta_{j,j}$ is a $(\alpha^{j+1}, \alpha )$-derivation and there exists $t_{2,j+1}$ such that 
$$
t_{2,j+1}at_{2,j+1}^{-1}=\alpha (a)+\delta_{n,j}(a)t_{2,j+1}^n+\ldots +\delta_{kn,j}(a)t_{2,j+1}^{kn}+\delta_{j+1,j+1}(a)t_{2,j+1}^{j+1}+\ldots  .
$$
%To prove that we can omit all the maps between $\delta_i$ and $\delta_{2i}$ in %the notations of Lemma, one can check that all these maps $\delta_j$, $n|j$ are %$(\alpha ,\alpha )$-derivations. Then one can find a parameter $t_2$ such that %new maps $\delta_j=0$. 

The rest of the proof is clear. For details see \cite{Zh, Prop.2, Cor.1, Lemmas 10, 3}. 
\qed
\endpf

\df Definition

Let $i=v_{K_2} (\varphi(t_2^n)(t_1)-t_1)\in n\Bbb N\cup \infty$,
($\varphi$ is defined in  subsection 8.1) 
and let $r\in \Bbb Z/i$ be  $v_{K_1}(x)\mo i$ 
where $x$ is  the residue of
$(\varphi(t_2^n)(t_1)-t_1)t_2^{-i}$.
Put   
$$
a=\res_{t_1}{\left(
\frac{(\delta'_{2i}-\frac{i+1}{2}{\delta'_i}^2)(t_1)}{{\delta'_i(t_1)}^2}dt_1 
\right )} \in K_0.
$$
($\delta'_i, \delta'_{2i}$ are the maps from the preceding lemma).
\enddf

\th Proposition

If $n=1$ then $i,r$ don't depend on the choice of a system of local parameters{{\rm;}}
if $i=1$ then $a$ does not depend on the choice of a system of local parameters{{\rm;}} 
if $n\not=1$ then $a$ depends only on the maps $\delta_{i+1}, \ldots
,\delta_{2i-1}$, $i,r$ depend only on the maps $\delta_j$, $j\notin n\Bbb N\,$, $j<i$. 
\endth

\pf Proof

We  comment on the statement first. The maps $\delta_j$ are uniquely defined by parameters $t_1, t_2$ and they depend on the choice of these parameters. So
the claim that $i,r$ depend only on the maps $\delta_j$, $j\notin n\Bbb N\,$, $j<i$ means that $i,r$ don't depend on the choice of parameters $t_1,t_2$ 
%such that they 
which 
preserve the maps $\delta_j$, $j\notin n\Bbb N\,$, $j<i$.  

Note that $r$ depends only on $i$. Hence it is sufficient to prove the proposition only for $i$ and $a$. Moreover it  suffices to prove it  for the case where $n\ne 1$, $i\ne 1$, because if $n=1$ then the sets $\{\delta_j: j\notin  n\Bbb N\,\}$ and $\{\delta_{i+1}: \ldots ,\delta_{2i-1}\}$ are empty. 

It is clear that $i$ depends on $\delta_j$, $j\notin n\Bbb N\,$.
Indeed, it is known that $\delta_1$ is an inner $(\alpha^2, \alpha )$-derivation (see the proof of the lemma). By \cite{Zh, Lemma 3} we can change a parameter $t_2$ such that 
$\delta_1$ can be made  equal $\delta_1(t_1)=t_1$. Then one can see that $i=1$. From the other hand  we can change a parameter $t_2$ such that 
 $\delta_1$ can be made equal to 0. In this case $i>1$. This means that $i$ depends on $\delta_1$. By \cite{Zh, Cor.3} any map $\delta_j$ is uniquely determined by the maps $\delta_q$, $q<j$ and by an element $\delta_j(t_1)$. Then using similar arguments and induction one deduces that $i$ depends on other maps $\delta_j$, $j\notin n\Bbb N\,$, $j<i$.   

Now we prove that $i$ does not depend on the choice of parameters $t_1,t_2$ 
%such that they 
which  preserve the maps $\delta_j$, $j\notin n\Bbb N\,$, $j<i$. 

Note that $i$ does not depend on the choice of $t_1$: indeed, if $t'_1= t_1 +bz^j$, $b\in K_1$ then $z^nt'_1z^{-n}= z^nt_1z^{-n}+(z^nbz^{-n})z^j = t'_1+r$, where $r\in \Cal M_K^i \backslash \Cal M_K^{i+1}$. One can see that the same is true for $t'_1=c_1t_1+c_2t_2^2+\ldots$, $c_j\in K_0$.  

Let $\delta_q$ be the first non-zero map for given $t_1,t_2$. If $q\ne i$ then by \cite{Zh, Lemma 8, (ii)} there exists a parameter $t'_1$ such that $zt'_1z^{-1}={t'_1}^{\alpha} +\delta_{q+1}(t'_1)z^{q+1}+\ldots$. Using this fact and Proposition 8.2 we can reduce the proof to the case where $q=i$, $\alpha (t_1)=\xi t_1$, $\alpha (\delta_i(t_1))=\xi \delta_i(t_1)$ (this case is equivalent to the case of $n=1$). Then we apply \cite{Zh, Lemma 3} to show that  
$$
\aligned 
v_{K_2} ((\phi (t'_2)-1)(t_1))&=v_{K_2}  ((\phi (t_2)-1)(t_1)), 
\endaligned
$$
for any parameters $t_2, t'_2$, i.e. $i$ does not depend on the choice of a parameter $t_2$. For details see \cite{Zh, Prop.6}.

To prove that $a$ depends only on $\delta_{i+1},\ldots ,\delta_{2i-1}$ we use the fact that for any pair of parameters $t'_1, t'_2$ we can find parameters $t''_1=t_1+r$, where $r\in \Cal M_K^i$, $t''_2$ such that corresponding maps $\delta_j$ are equal for all $j$. Then by \cite{Zh, Lemma 8} $a$ does not depend on $t''_1$ and by \cite{Zh, Lemma 3} $a$ depends on $t''_2=t_2+a_1t_2^2+\ldots$, $a_j\in K_1$ if and only if $a_1=\ldots =a_{i-1}$. Using direct calculations one can check that $a$ doesn't depend on $t''_2=a_0t_2$, $a_0\in K_1^*$. 

To prove the fact it is sufficient to prove it for $t''_1=t_1+ct_1^hz^{j}$ for any $j<i$, $c\in K_0$. Using \cite{Zh, Lemma 8} one can reduce the proof to the assertion that some identity holds. The identity is, in fact, some equation on  residue elements. One can check it by direct calculations. For details see \cite{Zh, Prop.7}. 
\qed
\endpf

\rk Remark

The numbers $i,r,a$ can be defined only for local skew fields
which splits. One can 
check that the definition can  not be extended to  the skew field in  Dubrovin's example. 
\endrk

\th Theorem

\Roster

\Item{(1)} $K$ is isomorphic to a two-dimensional local skew field
 $K_0((t_1))((t_2))$ such that 
$$t_2t_1t_2^{-1}=\xi t_1+xt_2^i+yt_2^{2i}$$
where $\xi$ is a primitive $n$th root, $x=ct_1^r$, $c\in K_0^*/(K_0^*)^{d}$,   $$y=(a+r(i+1)/2)t_1^{-1}x^2, \quad d=\gcdd(r-1,i).$$ 
If $n=1$, $i=\infty$, then $K$ is a  field.

\Item{(2)} Let $K,K'$ be two-dimensional local  skew fields of 
	characteristic zero which splits{{\rm;}} 
and let $K_1, K_1'$ be fields.
Let 	$\alpha^n=\id$, ${\alpha'}^{n'}=\id$ for some $n,n'\ge 1$.
	Then $K$ is isomorphic to $K'$
if and only if $K_0$ is isomorphic to $K_0'$ and the ordered sets $(n,\xi ,i,r,c,a)$ and $(n',\xi',i',r',c',a')$ coincide.
\endRoster
\endth

\pf Proof

(2) follows from the Proposition  of 8.2 and (1). We sketch the proof of (1). 

From Proposition 8.2 it follows that there exists $t_1$ such that $\alpha (t_1)=\xi t_1$; $\delta_i(t_1)$ can be represented as $ct_1^ra^i$. Hence there exists $t_2$ such that 
$$
t_2t_1t_2^{-1}=\xi t_1+xt_2^i+\delta_{2i}(t_1)t_2^{2i}+\ldots 
$$ 
Using \cite{Zh, Lemma 8} we can find a parameter $t'_1=t_1 \mo \Cal M_K$ such that 
$$
t_2t_1't_2^{-1}=\xi t_1+xt_2^i+yt_2^{2i}+\ldots
$$
The rest of the proof is similar to the proof of the lemma. Using \cite{Zh, Lemma 3} one can find a parameter $t'_2=t_2\mo \Cal M_K^2$ such that $\delta_j(t_1)=0$, $j>2i$. 
\qed\endpf

\th Corollary

Every two-dimensional local skew field $K$ with the ordered set $$(n,\xi ,i,r,c,a)$$ is a finite-dimensional
extension of a skew field with the ordered set  $(1,1,1,0,1,a)$.
\endth

\rk Remark

There is  a construction of a two-dimensional local skew field with a given
set $(n,\xi ,i,r,c,a)$.
\endrk

\eg Examples

\Roster

\Item{(1)} The ring of formal pseudo-differential equations is the skew field with the
set $(n=1, \xi=1, i=1, r=0, c=1, a=0)$.

\Item{(2)} The elements of $\Br(L)$ where $L$ is a two-dimensional local field of equal characteristic are local skew fields. If, for example, $L$ is a $C_2$- field, they split and $i=\infty$. Hence any division algebra in $\Br(L)$ is cyclic.
\endRoster
\endeg

\HHH 8.3.2. Characteristic $p$ case

\th Theorem

Suppose that a two-dimensional local skew field $K$ splits, $K_1$ is a field, $K_0\subset Z(K)$, 
$\chr (K)=\chr (K_0)=p>2$ and $\alpha=\id$.

Then $K$ is a finite dimensional vector space over its center if and only if $K$ is isomorphic to a two-dimensional local skew field
$K_0((t_1))((t_2))$ where
$$t_2^{-1}t_1t_2= t_1+xt_2^i$$ 
with $x\in K_1^p$, $(i,p)=1$.
\endth

\pf Proof 

The ``if'' part is obvious. We sketch the proof of the ``only if'' part. 

If $K$ is a finite dimensional vector space over its center then $K$ is a division algebra over a henselian field. In fact, the center of $K$ is a two-dimensional local field $k((u))((t))$. Then by \cite{JW, Prop.1.7} $K_1/(Z(K))_1$ is a purely inseparable extension. Hence there exists $t_1$ such that $t_1^{p^k}\in Z(K)$ for some $k\in \Bbb N\,$ and $K\simeq  K_0((t_1))((t_2))$ as a vector space with the relation 
$$
t_2t_1t_2^{-1}=t_1+\delta_i(t_1)t_2^i+\ldots 
$$ 
(see Remark 8.3.1). Then it is sufficient to show that $i$ is prime to $p$ and  there exist parameters $t_1\in K_1$, $t_2$ such that the maps $\delta_j$ satisfy the following property: 

\Roster

\Item{(*)}	
If  $j$ is not divisible by $i$ then $\delta_j=0$. 
If $j$ is divisible by $i$ then  
$\delta_j=c_{j/i}\delta_i^{j/i}$  with some $c_{j/i}\in K_1$.  
\endRoster 

Indeed, if this property holds then by induction one deduces that $c_{j/i}\in K_0$, $c_{j/i}=((i+1)\ldots (i(j/i-1)+1))/(j/i)!$. Then  one can find a parameter $t'_2=bt_2$, $b\in K_1$ such that $\delta'_j$ satisfies the same property and $\delta_i^2=0$. Then 
$$
{t'_2}^{-1}t_1t'_2=t_1-\delta'_i(t_1)t_2^i. 
$$

\smallskip

First we prove that $(i,p)=1$. To show it we prove that if $p|i$ then there exists a map $\delta_j$ such that $\delta_j(t_1^{p^k})\ne 0$. To find this map one can use \cite{Zh, Cor.1} to show that $\delta_{ip}(t_1^p)\ne 0$, $\delta_{ip^2}(t_1^{p^2})\ne 0$, $\ldots$, $\delta_{ip^k}(t_1^{p^k})\ne 0$. 

Then we prove that for some $t_2$  property (*) holds. To show it we prove that if  property (*) does not hold then there exists a map $\delta_j$ such that $\delta_j(t_1^{p^k})\ne 0$. To find this map we reduce the proof to the case of $i\equiv 1 \mo p$. Then we apply the following idea. 

Let $j\equiv 1 \mo p$ be the minimal positive integer
such that $\delta_j$  is not equal to zero on $K_1^{p^l}$. Then  one can prove that the maps $\delta_m$,  $kj\le m<(k+1)j$, $k\in \{1,\ldots ,p-1\}$ satisfy the following property:

\Roster
\Item{}
there exist elements $c_{m,k}\in K_1$ such that 
$$
(\delta_m-c_{m,1}\delta -\ldots -c_{m,k}\delta^k){|}_{K_1^{p^l}}=0
$$ 
where $\delta \colon K_1\rightarrow K_1$ is a linear map, $\delta {|}_{K_1^{p^l}}$ is a derivation, $\delta (t_1^j)=0$ for $j\notin p^l\Bbb N$, $\delta (t_1^{p^l})=1$, $c_{kj,k}={c(\delta_j(t_1^{p^l}))}^k$, $c\in K_0$.
\endRoster
\smallskip 

Now consider maps $\widetilde{\delta_q}$ which are defined by the following formula 
$$
t_2^{-1}at_2=a+\widetilde{\delta_i}(a)t_2^i+\widetilde{\delta_{i+1}}(a)t_2^{i+1}+\ldots , \quad a\in K_1. 
$$
 Then $\widetilde{\delta_q}+\delta_q+\sum_{k=1}^{q-1}\delta_k\widetilde{\delta_{q-k}}=0$  for any $q$.
In fact, $\widetilde{\delta_q}$ satisfy some identity which is similar to the identity in \cite{Zh, Cor.1}. Using that identity one can deduce that 

\noindent if  
\Roster
\Item{}  $j\equiv 1 \mo p$ and there exists 
the minimal $m$ ($m\in \Bbb Z$) such that  
$\delta_{mp+2i}{|}_{K_1^{p^l}}\ne 0$ if $j{\not |} (mp+2i)$
and 
$\delta_{mp+2i}{|}_{K_1^{p^l}}\ne s\delta_j^{(2i+mp)/j}{|}_{K_1^{p^l}}$ for any  $s\in K_1$ otherwise, and $\delta_q(t_1^{p^l})=0$ for $q<mp+2i$, $q\not\equiv  1\mo p$, 
\endRoster
\noindent then 
\Roster
\Item{} $(mp+2i)+(p-1)j$ is the minimal integer such that 
$\delta_{(mp+2i)+(p-1)j}{|}_{K_1^{p^{l+1}}}\ne 0$.
\endRoster

To complete the proof we use induction and \cite{Zh, Lemma 3} to show that there exist parameters $t_1\in K_1$, $t_2$ such that $\delta_q(t_1^{p^l})=0$ for $q\not
\equiv 1,2 \mo p$ and $\delta_j^2= 0$ on $K_1^{p^l}$.  
\qed\endpf

\th Corollary 1

If $K$ is a finite dimensional division algebra over its center then 
its index is equal to $p$. 
\endth

\th Corollary 2

Suppose that a two-dimensional local skew field $K$ splits, $K_1$ is a field, $K_0\subset Z(K)$, $\chr (K)=\chr (K_0)=p>2$, $K$ is a finite dimensional division algebra over its center of index $p^k$. 

Then either $K$ is a cyclic division algebra or has index $p$. 
\endth

\pf Proof

By \cite{JW, Prop. 1.7} $K_1/ \overline{Z(K)}$ is the compositum of a purely inseparable extension and a cyclic Galois extension. Then the canonical automorphism $\alpha$ has order $p^l$ for some $l\in \Bbb N$. By \cite{Zh, Lemma 10} (which is true also for $\chr (K)=p>0$), $K\simeq K_0((t_1))((t_2))$ with 
$$
t_2at_2^{-1}=\alpha (a) +\delta_i(a)t_2^i+\delta_{i+p^l}(a)t_2^{i+p^l}+\delta_{i+p^l}(a)t_2^{i+2p^l}+\ldots
$$   
where $i\in p^l\Bbb N$, $a\in K_1$. Suppose that $\alpha\ne 1$ and $K_1$ is not a cyclic extension of $\overline{Z(K)}$. Then there exists a field $F\subset K_1$, $F\not\subset Z(K)$ such that $\alpha {|}_F=1$. If $a\in F$ then for some $m$ the element $a^{p^m}$ belongs to a cyclic extension of the field $\overline{Z(K)}$, hence $\delta_j(a^{p^m})=0$ for all $j$. But we can apply the same arguments as in the proof of the preceding theorem to show that if $\delta_i\ne 0$ then there exists a map $\delta_j$ such that $\delta_j(a^{p^m})\ne 0$, a contradiction. 
We only need to apply \cite{Zh, Prop.2} instead of \cite{Zh, Cor.1} and note that $\alpha\delta =x\delta\alpha $ where $\delta$ is a derivation on $K_1$, $x\in K_1$, $x\equiv 1 \mo \Cal M_{K_1}$, because $\alpha (t_1)/t_1\equiv 1 \mo \Cal M_{K_1}$.

Hence $t_2at_2^{-1}=\alpha (a)$ and $K_1/ \overline{Z(K)}$ is a cyclic extension and $K$ is a cyclic division algebra $(K_1(t_2^{p^k})/Z(K), \alpha , t_2^{p^k})$. 
\qed\endpf

\th Corollary 3

Let $F=F_0((t_1))((t_2))$ be a two-dimensional local field, where $F_0$ is an algebraically closed field. Let $A$ be a division algebra over $F$.  

Then $A\simeq B\otimes C$, where $B$ is a cyclic division algebra 
of index prime to $p$ and $C$ is either cyclic {{\rm(}}as in Corollary 2{{\rm)}} or $C$ is a local skew field from the theorem of index~$p$. 
\endth

\pf Proof

Note that $F$ is a $C_2$-field. Then $A_1$ is a field, $A_1/F_1$ is the compositum of a purely inseparable extension and a cyclic Galois extension, and   $A_1= F_0((u))$ for some $u\in A_1$. Hence $A$ splits. So, $A$ is a splitting two-dimensional local skew field. 
 
It is easy to see that the index of $A$ is $|\Overline{A}:\Overline{F}|=p^qm$, $(m,p)=1$. Consider subalgebras $B=C_A(F_1)$, $C=C_A(F_2)$ where $F_1=F(u^{p^q})$, $F_2=F(u^m)$. Then by \cite{M, Th.1} $A\simeq B\otimes C$. 

The rest of the proof is clear. 
\qed\endpf

Now one can easily deduce that

\th Corollary 4

The following conjecture:
 the exponent of $A$ is equal to its index  for any division algebra $A$ over a $C_2$-field $F$ {{\rm(}}see for example \cite{PY, 3.4.5.}{{\rm)}} 

\noindent has the positive answer for $F=F_0((t_1))((t_2))$. 
\endth    

\vskip 1cm

\Bib   Reference 

\rf {JW} B. Jacob  and  A. Wadsworth, Division algebras over Henselian fields, J.Algebra,128(1990), p. 126--179.  

\rf {M} P. Morandi, Henselisation of a valued division algebra, J. Algebra, 122(1989), p.232--243.

\rf {Pa} A.N. Parshin, On a ring of formal pseudo-differential operators, Proc. Steklov Math. Inst., v.224 , 1999, pp. 266--280, (alg-geom/ 9911098).  

\rf {PY} V.P. Platonov and V.I. Yanchevskii, Finite dimensional division algebras, VINITY, \linebreak 77(1991), p.144-262 (in Russian)

\rf {Zh} A. B. Zheglov, On the structure of two-dimensional local skew fields,
to appear in Izv. RAN Math. (2000). 

%\rf {Ma} H. Matsumura, Commutative ring theory, Cambridge Univ. Press, 1992. 

\endBib

\Coordinates

Department of Algebra, Steklov Institute,
Ul. Gubkina, 8, Moscow GSP-1, 117966 Russia.

E-mail: abzv24\@mail.ru, azheglov\@chat.ru
\endCoordinates

\vfill
\pagebreak

\end